%% file: relkren.plain.tex
\input macros.tex

\input amssym.tex
\input amssym.def
\magnification=1400
\baselineskip = 15pt
{\centerline 
{\bf 
Partitions with 
independent iterates in random dynamical systems}}
 \bigskip
\centerline{Boris Begun (1)}

\medskip\centerline{Andr\'es del Junco (2) }
\vskip1in

\centerline{\sl Abstract}
Consider an invertible measure-preserving transformation of
a probability space.
A finite partition
of the space is called
weakly independent if there are
infinitely many images of this partition
under powers of the transformation
that
are jointly independent.
Krengel proved that 
a transformation is weakly mixing if and only if
weakly independent partitions 
of the underlying space
are dense among
all finite partitions.
Using the tools developed in the later
papers of 
del Junco-Reinhold-Weiss and 
del Junco-Begun we obtain 
Krengel-type results for weakly mixing
random dynamical systems 
(or equivalently, skew products that are relatively
weakly mixing).
 \vskip1in
\item{(1)} Dept. of Mathematics, Hebrew University  \hfill\break
             begun@math.huji.ac.il
\medskip
\item{(2)} Dept. of Mathematics,  University of Toronto\hfill \break
              deljunco@math.toronto.edu
\vfill\eject

\noindent {\sl 0. Introduction}

The notion of weak mixing is one of the strengthenings
of the notion of ergodicity of a measure-preserving transformation.
Along with the notions of mixing, $K$-property etc. it belongs 
to the hierarchy
of statistical properties of transformations in  ergodic theory.
This notion complements the notion of discrete spectrum and
has several equivalent characterizations mostly going back to Koopman and
von Neumann (early 1930's).
In  1970 Krengel proved in [K] another
unexpected
characterization
of weak mixing.
In this paper we prove two analogues of Krengel's result for systems that are
weakly mixing relative to a factor.
Our claim can be alternatively stated in the
language of random dynamical systems.
Before presenting the original
theorem and its relativization, we provide necessary
definitions concerning partitions of probability spaces
(a partition is a representation of the space as a finite disjoint union
of measurable sets).

A family of measurable sets in a probability space 
($X,m$) is independent 
if for every finite sub-family $A_1,A_2, \ldots, A_k$ we have
$$
m(\bigcap_{i=1}^{k}A_i)=\prod_{i=1}^{k} m(A_i).
$$
A family of measurable partitions in $X$ is independent 
if  every family of sets of the partitions -- one set from each 
partition -- is independent.
If $P=\{A_1, \ldots, A_n\}$ is an ordered partition then 
${\rm dist}\, P$,
the distribution of $P$, 
is a vector in $ \R^n$ defined by $({\rm dist}\, P)_i=m(A_i)$.
Consequently, 
if $P=\{A_1, \ldots, A_k\}$ 
and $Q=\{B_1, \ldots, B_k\}$ are two partitions 
with the same number of atoms, we will say that $P$ and $Q$ 
have the same distribution if $m(A_i)=m(B_i)$ for all $i$.
The set of all (ordered) partitions into $n$ atoms can be turned
into a complete metric space with  
the metric  
$d_m(P,Q)=\sum_{i=1}^{n}m(A_i\triangle B_i)$.
Wherever there is no ambiguity we will omit the subscript and denote
this distance by plain $d$.

\def\sX{(X, \c B, m, T)}
 
By a (dynamical) system we mean a quadruple
${\bf X} = \sX$ where $(X, \c B, \mu)$ is a probability space and
$T$ is an sutomorphism of  $(X, \c B, \mu)$, that is $T$
  is an invertible $\mu$-preserving map.
A sequence $\set {n_i} \subset \Z$ is called mixing for $T$ if
$$
\mu(A \cap T^{n_i}B) \to \mu(A) \mu(B) \hbox{ as } i \to \infty.
$$
One characterization of   weak  mixing of $T$ is that the
product system $(X \times X, \c B \times \c B, \mu \times \mu, T \times T)$
is ergodic. An equivalent one is that there exists a mixing sequence
$\set{n_i} \subset \N$ for $T$ having density one in $\N$ or equivalently
there exists any mixing sequence at all.

\proclaim Theorem (Krengel).
$(X, \c B, m, T)$ is a weakly mixing system if and only if for every 
finite measurable partition $P$ of $X$ and $\ep >0$,
there is an infinite sequence $0=n_0<n_1<\ldots$ and
a partition $P'$ with the same distribution as $P$ such that
$d(P,P') < \ep$ and the partitions 
$\{T^{-n_i}P'\}_{i=0}^{\infty}$
are jointly independent.

 Krengel's  very technical proof was 
substantially simplified in 1999 by del Junco, 
Reinhold and Weiss ([JRW]), paving the way 
for generalizations and refinements.
Among other things it was shown there that the sequence in Krengel's theorem
can be chosen inside any prescribed mixing sequence for $T$.
The key tool for this was an extension result for stationary measures
on $A^{\Z}$.
In  the subsequent paper [BJ]  the extension theorem was generalized to
$A^G$, $G$ any discrete amenable group, permitting the Krengel theorem 
to be generalized to
free actions of discrete amenable groups.

\def\sX{(X,\c B,\mu,T)}

\def\ph{\phi}

\def\prP{ \Bbb P}
\def\Om{\Omega}
\def\om{\omega}

We restrict our attention to systems $\sX$
 such that $X$ is a complete metric space and
 $\c B$ is the $\mu$-completion
 of the
Borel $\si$-algebra of $X$.  We will call such systems standard.
This assumption 
involves no loss of generality for theorems of the
kind we are concerned with here.
We will be studying the properties of $\X$ relative to a
   distinguished factor algebra ${\c F}$.  There is no loss of generality
in assuming that $\c  F$ is given by
 ${\cal F} = \pi^{-1} \bar{\c F} $,
where
$$
\pi: {\bf X} \to {\bf \Om} = {(\Omega, \bar {\c F}, \prP ,\theta)}
$$
is a factor map from ${\bf X}$ to another standard system
${\bf \Omega}$.
 The sets $\pi^{-1}\set{\om}, \ \om \in \Om$, are the fibers of 
$\c F$.
 The measure $\mu$ decomposes over the
 factor as $$\mu = \int_\Om \mu_\om d \prP (\om),$$
 where $\mu_\om$ is supported on the fiber $\pi^{-1}\set{\om}$.

As is well-known, when  the factor system  is ergodic, 
we can represent $\bf X$ as a skew-product over $\bf \Om$:
$$
\eqalign{
(X,\c B,\mu)=& (\Omega, \bar {\c F}, \prP) \times (Y, \c C, \nu),   \cr
 {\c F} =& \bar {\c F} \times Y,   \cr
T(\omega, y)=& (\theta\omega,T_{\omega}y),
}
$$
where $T_\om$ is an automorphism of $ (Y, \c C, \nu)$ which depends
measurably on $\om$.
The pair $(\th, \set{T_\om})$ is called a random dynamical system.

 From now on we identify $\c F$ with $\bar{\c F}$
and write $\c F$ for both.
One definition of relative weak mixing uses the relative product of
$\sX ={\bf X}$ with itself. This is the 
system
 $$
{\bf X} \times_{\bf \Omega}{\bf X} =
(X \times X, \c B \times \c B,
\lambda,
T \times T)
$$ 
where the relative product measure 
$\lambda = \mu\times_{\bf \Om} \mu$ is defined by 
$$\lambda = \int_\Om \mu_\om \times \mu_\om d \Bbb P(\om).$$
Note that $\la$ is supported on the fibered product
$$
\cup\set{\pi^{-1}\set{\om} \times \pi^{-1}\set{\om}: \om \in \Om} \subset X \times X.
$$
Consequently ${\bf \Om}$ is a factor of ${\bf X} \times_{\bf \Omega}{\bf X}$
in a canonical way.

 By definition
${\bf X}$  is  weakly mixing relative to  $\c F$
if ${\bf X} \times_{\bf \Omega}{\bf X}$
 is
   ergodic relative to $\c F$.
What this means is that any invariant function for the relative product is
$\c F$-measurable. 
  There is an equivalent characterization as a (relatively)
 mixing condition: $T$ has
a relatively mixing sequence $\set{n_i}$,  which means
 that for any sets
$A$ and $B$
$$
\mu_\om(A \cap T^{-n_i}B) - \mu_\om(A) \mu_\om(T^{-n_i}B)
$$
converges to $0$ in measure (as a function of $\om$).
 Equivalently, $T$ has a relatively mixing sequence of density one in $\N$.
The existence of a relatively mixing sequence of density one in the presence
of relative weak mixing follows easily
from [F], Lemma 7.6. The converse fact, that a relatively mixing sequence
implies relative weak mixing, is easy since a relatively mixing sequence
for ${\bf X}$ is also a relatively mixing sequence for
 ${\bf X} \times_{\bf \Om} {\bf X}$ and the existence of a relatively
mixing sequence certainly forces relative ergodicity.
  
We will present two versions of Krengel's theorem for relatively 
weakly mixing systems.
The first says that by a uniformly small 
perturbation of a  given partition it is possible to obtain 
one with the same distribution on  a.a. fibers so
that it has infinitely many independent iterates  on a set of
fibers with probability 
arbitrarily close to 1. ``Probability''
here means the measure on the factor, in accordance with the 
ideology of RDS theory.
Consequently, ``a.s.'' will mean 
``for $\prP$-a.e. $\omega \in \Omega$''.
 
 For convenience we will state and prove 
the results  in the case of ergodic systems.
Theorem 1 remains valid without this assumption and so does Theorem 2, with
a minor modification. The proofs are not much more
difficult. After proving the results in section 3 we will make some comments
about how to remove the ergodicity assumption. We remark that when
$T$ is weakly mixing relative to $\Om$, ergodicity of $T$ is equivalent to
ergodicity of $\th$.

\proclaim Theorem 1.
Let 
$(X,\cal B,\mu,T)$ be an ergodic  system with a
factor 
$(\Omega, \cal F, \prP ,\theta)$ 
where $\mu$ decomposes as $\mu=\int \mu_\omega d\prP$.
Assume $\{n_i\}_{i=0}^{\infty}$ is a relatively mixing sequence for
$T$, and $P$ is a finite 
partition of $X$ for which there is an $\al >0$ such that
$$
\forall p \in P:
 \mu_{\omega}(p)>\alpha >0
{\ \rm for\ }
 \prP {\rm-a.e.\  } \omega\in \Omega. 
$$
Then for every $\ep>0$ there exist a partition
$Q$ of $X$,
a set $E \subset \Omega$, $\prP (E) < \ep$,
and 
a subsequence $\{m_i\}$ of $\{n_i\}$, such that
$$
\dist_{\mu_{\omega}} Q = \dist_{\mu_{\omega}} P 
{\hbox { for }} \prP{\hbox{-a.e. }}\omega 
\eqno (1)
$$
$$
d_{\mu_{\omega}} (P,Q) < \ep \hbox { for } \prP \hbox{-a.e. }
\omega
\eqno (2)
$$
$$
\set{T^{m_i}Q  }\hbox{ is jointly independent  }
\hbox {  with respect to }
\mu_{\omega}\  \forall \omega \not \in E
\eqno (3)
$$

\def\P{\Bbb P}
 
 The independence claim holds for a set of $\omega$ of 
measure  arbitrarily close to 1, not $\prP$-almost everywhere.
This is not a drawback of the method --
if we require that $Q$ be uniformly close to $P$ then we cannot
obtain independence on almost all fibers.
For a counterexample consider the so-called
$(S,S^{-1})$ transformation in its simplest version.
Let $\Omega=Y=\set{-1,1}^{\Z}$ with the standard product measure $\P$
coming from $(1/2, 1/2)$ distribution of probabilities on 
$\{-1,1\}$, and let $\th=S$ be the Bernoulli shift: $(Sx)_n=x_{n+1}$.
The $(S,S^{-1})$ map is a skew product $T$ on $\Om \times Y$ defined by
 $$T(\omega,y) = (S\omega, S^{\omega_0}y).$$ The 
family $\set{S^{\omega_0}}$ 
contains only two
distinct transformations -- the shifts to the left and to the right.
The system is weakly mixing with respect to the 
base (it is even strongly mixing in the $L^1$ sense -- see [R] for the
definitions and the proof), but
the independence a.e. cannot be achieved even 
along a set of two numbers, say $\set{n_1,n_2}$.

To see this first observe that if $T^{n_1}Q$ were
$\mu_\om$-independent of
$T^{n_2}Q$  for a.a. $\om$ then, for $n = n_1-n_2$, we would also have 
 $T^{n}Q$
$\mu_\om$-independent of
$Q$ for a.a. $\om$. Hence we may as well take $n_1 = 0$ and $n_2 = n$.
Choose a partition
 $R$ of $Y = \Om$ such that $$\dist_{\P} R = (1/2,1/2) \hbox{ and }
 d_\P (R,SR) < {1 \over 100}.$$
 Example: viewing $R$ as a map into the index set $\set{-1,1} $  let
$R(y)={\rm sign}(\sum_{i=1}^{i=10^6 + 1}y_i).$
 (Note that this sum cannot be $0$.)
Define $ P$ on $X$ by $P(\om,y) = R(y)$.
Observe that 
$$T^n(\om,y) = (\th^n \om, S^{\ph(n,\om)}y)$$ where
 $$\ph(n, \om) = \sum_{i=0}^{n-1}\om_i.$$
Now let $\de = 0$ or $1$ according to whether $n$ is odd or even.
Then the set $E = \set{\om:\ph(\om,n)= \de}$  has positive measure
and for $\om \in \th^n E$ we will have 
$$d_{\mu_\om}( P,T^{ n}  P) = d_{\P}(R,S^\de R) < {1 \over 100}.$$
Now suppose $Q$ is a partition of $X$ such that for almost all $\om$ 
$$\dist_{\mu_\om}(Q) = ({1 \over 2},{1\over 2})\hbox{ and }
 d_{\mu_\om}(Q,P) < {1 \over 100}.$$ Then it is easy to see that 
for $\om \in \th^n E$
$$d_{\mu_\om}( Q,T^{ n}  Q) < {3 \over 100}.$$
On the other hand if $Q$ and $T^nQ$ were independent with respect to $\mu_\om$ then we would have
$d_{\mu_\om}(Q,T^nQ) = {1 \over 4}$, a contradiction.
\bigskip
In our second result, by dropping  the requirement that
the perturbation be uniformly small we are able to obtain the
desired independence on almost all fibers, rather than just a large
set of fibers.

\proclaim Theorem 2.
With the same hypotheses as in Theorem 1
  one can find 
a partition $Q$  and a subsequence $\{m_i\}$ of $\set{n_i}$ such that
$\dist_{\mu_\om} Q = \dist_{\mu_\om} P$ for a.a. $\om$,
$d_\mu(P,Q) < \ep$ and the partitions
$\set{T^{m_i}Q}$ are jointly independent with respect to $\mu_\om$ for
almost every $\om$.

\noindent{\sl Remark.}\ Of course to say that $d_\mu(P,Q)$ is small is just
to say that $d_{\mu_\om}(P,Q)$ is small for all but a small set of $\om$.

\medskip

Our presentation conforms to the following plan. 
After a section on preliminaries and notation we prove
in Section 2 a result (Propostion 1) on the existence of measures with 
prescribed marginals -- a non-stationary generalization of the stationary
 extension result Theorem 2 in [BJ]. It's proof mimics closely the proof of
the stationary result. Theorems 1 and 2 are proved in Section 3.

\bigskip

\noindent {\sl 1. Notation and Preliminaries}

\bigskip
In this section we introduce some notation and definitions
that pertain to Sections 2 and 3. In the beginning of
Section 3 more preliminaries are collected that are specific for
the proof of the main result.

First we adopt a more formal definition of an ordered partition.
A (finite) partition of a measurable space $(X,\cal B)$ is a measurable 
map $P$ from $X$ to a finite 
index set $A$ (so by definition 
we are dealing with ordered partitions). The sets $p=P^{-1}(a)$, 
$a \in A$,
are the atoms of $P$ and we write $p \in P$. If  $(X,\cal B)$  
carries a probability measure
$\mu$ then $\dist P$ 
(or $\dist_{\mu}P$ if the measure needs to be emphasized)
denotes the measure $\mu \circ P^{-1}$
on $A$. If $B \in \cal B $, $\mu (B) \ne 0$ 
then $\dist_{\mu}(P|B)$ refers to the restriction
$P|_{B}$ and 
the normalized 
measure $\mu_{B}$. We write $P \prec Q$ ($Q$ refines $P$)
if each atom of $P$ is a union of atoms of $Q$.
For two partitions $P$ and $Q$ indexed by the same alphabet $A$ the 
distance $d(P,Q)$, or $d_{\mu}(P,Q)$, between them is 
$$d(P,Q)= \mu\{x \in X \mid P(x) \ne Q(x)\}.$$
For a fixed alphabet $A$ this metric defines 
a complete metric space of partitions.

Two partitions of a  probability
space $(X,\mu)$ are independent if every atom $p$ of the 
first one is independent of every atom $q$ of the second
one: $\mu (p \cap q) = \mu(p) \mu (q)$.
The definition of joint independence of more than two partitions
is similar.
If $P$ and $Q$ are finite partitions of a probability 
space $(X,\mu)$ then we will say 
$P$ is $ \delta$-independent of $Q$
whenever
 $$
\|\dist (P|q) - \dist P\|_{\infty} < \delta\  \forall
q \in Q.
$$
(We are identifying a measure $\si$ on
a finite set $A$ with the vector $\{\si(a)\}_{a \in A}$.)
When $P$ is $\delta$-independent of $Q$ we will 
write $P \perp_{\delta} Q$
(or $P \perp_{\delta , \mu} Q$).
Note that this definition of approximate independence
differs from a more standard one which uses the $l_{1}$-norm
and allows a small exceptional set of $q$'s.
In fact any reasonable defintion of approximate independence could be
used in the statements of our results but the one we have chosen is convenient
for the proofs.
We will make use of the fact that if  $P \perp_{\delta} Q$ and
$Q \succ R$  then $P \perp_{\delta} R$.

One more piece of notation: we write $x \close \de y$ whenever $x$ and $y$
are two elements of a space which are no more than $\de$ apart with respect
to some metric or norm which is clear from the context.
\bigskip


\noindent {\sl 2. Building a measure with prescribed marginals}

\bigskip

This section is static -- no dynamics are involved. The 
object of study is measures on $A^{\Sigma}$, where $A$ is a finite
set (alphabet) and 
$\Sigma$ is countable. We prove a result on the existence of a
probability measure with prescribed projections
 on $A^K$ for $K \in \c F$, a certain family of finite subsets 
$\cal \F$ of $\Sigma$.

We will be working with
measures on $A^H$ for various subsets $H \subset \Sigma$ and we use the
usual Borel structure on $A^H$. Mostly $H$ will be finite.
Let $K_1, K_2 \subset \Sigma$.
Whenever $K \subset L\subset \Sigma$ we will 
denote by $P^{K}$ the projection map from 
$A^{L}$ to $A^{K}$, which is a finite partition of $A^{L}$ when $K$ is
finite. We will also write $P^K = \pi_K$. $\pi_K$ also acts on measures:
if $m$ is a measure on $A^L$ then $\pi_K m = m \circ \pi^{-1}$.  
$m$ is then an extension of 
$\pi_K m$. 
By definition, if $E \subset A^K$, then $m(E)$ means the same as $\pi_K m (E)$.
 Two measures on $A^{K_1}$ and $A^{K_2}$ are
called consistent if their projections on $A^{K_1 \cap K_2}$ coincide.
Consistency of a family of measures means pairwise consistency.

For $H$ finite we will say a probability 
measure $m$ on $A^{H}$ 
is {$\delta$-independent} whenever there is an indexing 
$H=\{h_{1}, \ldots, h_{n}\}$ such that 
$$P^{\{h_{i}\}} \perp_{\delta,m} P^{\{h_{1}, \ldots, h_{i-1}\}}\ 
 \forall i=2, \ldots , n.$$
One easily sees that $\delta$-independence 
of $m$ implies 
$$
\|m - p_H\|_{\infty} < (|H|-1)\delta,
$$
where $p_H$ is the product measure on $A^{H}$ that has the same  
one-dimensional marginals as $m$. 
Note also that if $m$ is $\delta$-independent then so is 
$\pi_I m$ for any $I \subset H$.


\proclaim Proposition 1.
For every $N \in \N$ and $\alpha>0$ there is a $\delta=\delta(\alpha,N)>0$ 
such that 
the following assumptions imply the conclusion below.
\vskip.1in
\noindent Assumptions:
\vskip.01in
$A$ is a finite alphabet, $\cal F$ is a family of finite 
subsets of $\Sigma$ and  for every set $K\in \cal F$ a probability 
measure $\mu_K$ on $A^K$ is given. These data satisfy:
\vskip .05in
\itemitem{ (F)}
$
\forall n \in  \Sigma \hbox{ we have } |\cup\{K\in \c F : n \in K\}| \le N
$
\itemitem{(M1)} 
The measures $\mu_K$ are pairwise consistent, 
\itemitem{(M2)} 
Every atom of each 1-dimensional marginal of each measure $\mu_K$  
has measure  $\ge \alpha$.
\itemitem{(M3)} 
Every measure $\mu_K$ is $\delta$-independent.
\vskip.05in
\noindent Conclusion:
\vskip.05in
 The measures 
$\{\mu_K\}_{K\in\cal F}$ 
have a common extension 
to $A^\Sigma$.



\noindent {\sl Remark}.\quad
Evidently there would be no loss of generality in  assuming  in the proposition that
$\cup \c F = \Sigma$.

 Proposition 1 is an immediate consequence  of  the following
 claim. 

\proclaim Claim A.
For every $N \in \N$ and $\alpha>0$ there is a $\delta=\delta(\alpha,N)>0$ 
and a $\beta=\beta(\alpha,N)> 0$ such that 
the following implication holds.
Suppose that the assumptions of Theorem 1 hold.
Suppose further that $I$ is a finite subset of $\Sigma$ and
$\lambda$ is
a measure on $A^{I}$ that is consistent with every $\mu_K$ 
and $\beta$-independent.
Then there is a 
measure $\lambda'$ on $I'$ that extends $\lambda$, 
and is again consistent with every $\mu_K$ and $\beta$-independent.

\noindent {\sl Remark.} Clearly $\delta \le \beta$.
Typically $\delta \ll \beta \ll 1$.
 Explicitly, one can 
take $\beta=\alpha^N/2N$, as will be seen from the proof. 
The expression 
for $\delta=\delta(\alpha,N)$ can also be developed in principle
but we will not need it.

\proof
To prove  Claim A we will need two lemmas. These are the 
Lemmas 2.1 and 2.2 from [BJ], the first of which appeared
already 
in [JRW] (Proposition 1.2). To make our presentation 
self-contained, we provide \break sketches of the  proofs 
here and refer the reader to the previous papers for more 
details.

The first lemma deals with {\sl signed} measures, for which 
the notions of restriction, extension and consistency remain meaningful.

\proclaim Lemma 1.
Let $K$ be a finite set and $\{K_i\}$  a finite family of subsets of $K$ 
whose union is $K$.
Any consistent family of signed measures on $A^{K_i}$ has a 
common extension -- a signed measure on $A^{K}$.

\def\J{\c J}
\noindent{\sl Sketch of the proof.}
We are given the measures $\mu_{K_i}$, $i=1,\dots,n$. 
Pick an arbitrary {\sl probability} measure $q$ on $A^K$.
For every 
$\J\subset\{1,\dots,n\}$ denote by $\mu_{\J}$ the common projection of 
$\mu_{K_j}$, $j \in \J$, onto 
$A^{\cap\set{K_j: j \in \J}}$, multiplied 
by $\pi_{\set{K_j: j \in \J}}q$.
One common extension is
$$
\mu = \sum_{\varnothing \ne {\cal J} \subset \{1,\dots,n\}}{(-1)^{|{\J}|+1} \mu_{\J}}.
$$
This expression is inspired by the inclusion-exclusion 
formula in combinatorics and the proof is based on 
the same idea -- when
the expression for $\mu$ is projected to any of 
the $A^{K_i}$'s the sum over  all
${\cal J} \neq \set{1, \ldots n}$
 can be divided 
into pairs that cancel out.

\endproof

\proclaim Lemma 2.
Suppose $V$ and $W$ are finite-dimensional real normed vector
spaces and $\Pi \colon V \to W$ is any surjective linear map. 
Then there is a
constant $C=C(\Pi)$ with the following property: if $v \in V$ 
and
$w=\Pi(v) \neq 0$ then $\Pi$ has a right inverse $B$
($\Pi B =id_{W}$) such that $Bw=v$ and
$\|B\| \leq  C 
{\|v\|}/ {\|w\|}$.

\noindent {\sl Sketch of the proof.}
Using the equivalence of any two norms in a finite-dimensional 
space it is not hard to see that there is no limitation 
of generality in the assumption that $V$ is Euclidean 
and $\Pi$ is the orthogonal projection onto its 
subspace $W$. In this  case the result is immediate ($C=1$). 

\endproof

Continuing with the proof of 
Claim A, 
we shall assume, as we may, that $\bigcup \c F = \Sigma$.
Fix $n \in \Sigma \backslash I$, denote $I'=I\cup\{n\}$, 
${\cal F}_n=\{K\in{\cal F} \mid n \in K\}$ and consider the
following two families 
of subsets of $\Sigma$:
$$
{\cal S}=\{K\cap I' \mid K \in \c F_n \},
$$
$$
\c R=\{K\cap I \mid K \in \c F_n \}=
\{S\backslash \{n\} \mid S \in \c S\}.
$$
$\c S$
and $\c R$ are finite families of sets.
All sets in $\c S$ contain $n$.
Let $$\bar S=\cup \c S \hbox{ and }
 \bar R=\cup\c R.$$
Clearly $\bar R=\bar S\backslash \{n\}$, $\bar S \subset I'$, 
$\bar R \subset I$, and by $({\sl F}
)$
the cardinality of $\bar S$ does not exceed $N$.

We are going first to construct a probability
measure $\si$ on $A^{\bar S}$ -- a prospective marginal of $\la'$.
It necessarily has to satisfy the three conditions derived from
the corresponding conditions imposed on $\la'$ by Claim A, namely

$$
 \pi_{\bar R}\si=\pi_{\bar R}\lambda ,  \eqno  (1)
$$

$$
\pi_{K\cap I'}\si=\pi_{K\cap I'}\mu_{K}\qquad \hbox 
{for every $K\in \c F_n$}, \eqno  (2)
$$

$$ 
 P^{\{n\}} \perp_{\beta,\si} P^{\bar R}. \eqno  (3)
$$

Let us assume for the moment that we have already constructed 
$\si$ satisfying
(1), (2) and (3).   Because of (1)
we can then 
define $\la'$ on $A^{I'}$ to be the relative product
measure
$\la' = \la \times_{ A^{\bar R}} \si$, that is

$$
\la'(x) = {\la(\pi_I x) \si(\pi_{\bar S} x) \over
\la(\pi_{\bar R} x)} \qquad \forall x \in A^{I'}.
$$
$\la'$ is certainly an extension of $\la$ and of $\si$.
$\la'$ is consistent with every $\mu_K$, $K \in \c F$,
or equivalently 
$\la'$ is an extension of $\pi_{K\cap I'}\mu_{K}$,
because  either $\si$ or $\la$ is such an extension
(depending on whether $n \in K$ or not).

Moreover the definition of $\la'$ 
together with (3) implies that for all $X \in P^{\{n\}}$, 
$Y \in P^{\bar R}$ and
$Z \in P^{I \backslash \bar R}$ 
we have
$$
\la'(X|Y \cap Z)=\si(X|Y) 
\close {\beta} 
\si(X)=\la'(X).
$$ 
This means 
that $P^{\{n\}} \perp_{\beta,\la'}P^{I}$, 
and since
$\la$ is $\beta$-independent it follows that 
$\la'$ is also $\beta$-independent,
completing the proof 
of Claim A and  Proposition 1.

We now proceed to construct $\si$. 
For each $S = K \cap I' \in \c S$ and for each
$R = K \cap I \in  \c R$, where $K \in \c F_n$,
 let 
$$
\mu_S=\pi_{S}\mu_{K} \hbox{ and }
\nu_{R}=\pi_{R}\mu_{K}.
$$
By the consistency of the measures $\mu_K$ these 
requirements 
unambiguously define $\mu_S$ and $\nu_R$ even though $K_1 \cap I' = K_2 \cap I'$ does not imply 
$K_1 = K_2 $. Moreover,
both $\{\mu_S\}_{S\in \c S}$ and
$\{\nu_{R}\}_{R\in \c R}$ 
are consistent
families of measures.
We let $V$ denote 
the vector
space of all signed measures $\rho$ on $A^{\bar R}$ and $W $ the space 
of all consistent families
 $\{\rho_{R}\}_{R \in {\c R}}$, where
 $\rho_{R}$ is a signed 
measure on $A^R$, 
both  spaces 
endowed  with the $l_{\infty}$-norm.
 We let $\Pi: V \to W$ 
denote the projection map, that is 
$\Pi(\rho)= \{\pi_R\rho\}_{R \in \c R}$ 
for each  $\rho \in V$, so we have
$\Pi (\pi_{\bar R}\lambda)= \{\nu_{R}\}_{R \in \c R}$. 
Lemma 1 tells us that
 $\Pi$ is surjective
 so by Lemma 2 we conclude that
$\Pi$ has a right inverse 
$B$ 
such that 

$$
B(\{\nu_R\}_{R \in \c R})=\pi_{\bar R}\lambda,$$
 and 

$$
\| B \|\le C\| {\pi_{\bar R}\lambda}\|/\|{\{\nu_R\}_{R \in \c R}}\|,
$$
 where
$C=C({\Pi})$ is a constant depending  on $\Pi$.
However, 
$\Pi$ is completely determined by 
specifying the sets $R \in \c R$, which are
subsets of $\bar R$, and by the condition ($\c F$)
the cardinality of the latter is less than $N$, so 
there are less than $2^N$ subsets.
Therefore there are less 
than $(2^N)^{2^N}$
possibilities for 
the projection $\Pi$. This means that we may 
take $C$ to be a constant depending only on $N$.
Since $\pi_{\bar R}\la$ is a probability measure we have 
$\|{\pi_{\bar R}\la\|}_{\infty} \leq 1$. Moreover the 
hypothesis (M2) ensures that 
 $\abs A \leq  \alpha^{-1}$.
Since each $\nu_{K}$ is a 
probability measure on $A^K$, a 
set of cardinality less than $\abs A^n$, we have
$$
\|{\{\nu_R\}_{R \in \c R}}\|_{\infty} > 
\alpha^{N},
$$
so we obtain an absolute bound
$\|{B}\| \leq C\alpha^{-N}= C'$.

Recalling that $\bar R \cup \set n = \bar S$,
we view each $y \in A^{\bar S}$ 
as a pair $y=(x,a),\  x \in A^{\bar R},\ a \in A$,
and use a similar convention for the sets $A^S$, $S\in \c S$. 
With this convention,  if $\si$ is a signed
measure on $A^{\bar S}$ and $a \in A$ then $\si(\cdot,a)$ is a (signed)
measure on $A^{\bar R}$.
Note that the family $\{\mu_{S}\}_{S\in \c S}$  
can be rewritten in the form
$\{\mu_{R\cup\{n\}}\}_{R\in \c R}$,
and $\mu_{R\cup \{n\}}(\cdot,a)$ is a measure 
on $A^R$ for each $R \in \c R$.
We now define
$\si$ by specifying that
$$
\si(\cdot,a)= B(\{\mu_{R\cup\{n\}}(\cdot,a)\}_{R \in \c R}) 
\quad  \forall a \in A.
$$
Note that the consistency of the family 
$\{\mu_{R\cup\{n\}}(\cdot,a)\}_{R \in \c R}$
follows from the consistency of the family 
$\{\mu_{R\cup\{n\}}\}_{R \in \c R}$.
We then have
$$
\eqalign{
\pi_{\bar R} \si
 &=\sum_{a \in A} \si (\cdot,a)  \cr
 &=\sum_{a \in A} B (\{\mu_{R\cup \{n\}}(\cdot,a)\}_{R \in \c R}) \cr
& = B\left(\sum_{a \in A} \{\mu_{R\cup\{n\}}(\cdot,a)\}_
{R \in \c R}\right)\hbox{ \quad\quad (since {\sl B} is linear)} \cr
& = B(\{\nu_{R}\}_{R \in \c R}) = \pi_{\bar R} \la,
}
$$
establishing (1). To check (2) observe that for each 
$a \in A$
and $R \in \c R$ we have
$$
(\pi_{R\cup\{n\}}\si)(\cdot,a)= \pi_{R}(\si(\cdot,a))
=\pi_{R}B(\{\mu_{Q\cup \{n\}}(\cdot,a)\}_{Q \in \c R}) = \mu_{R}(\cdot,a),
$$
since $\Pi B={\rm id}_W$.

It remains to check that $\si$ is non-negative and satisfies (3).
Note that all $\nu_{R}$ and $\mu_{R\cup \{n\}}$, $R \in \c R$,
are $\delta$-independent. 
Denote generically by
$p_L$ the product measure of 1-dimensional
marginals on $A^{\{l\}}, l \in L$. In particular $p_{\{n\}}$ is the common
1-dimensional marginal of all measures in $\c F_n$, ).
Now recall that for $R \in \c R$

$$
\|\mu_{R\cup \{n\}}-p_{R\cup \{n\}}\|_{\infty} < 
(|R\cup \{n\}|-1)\delta =|R|\delta< N\delta.
$$

\noindent Similarly $\|\nu_R-p_R\|_{\infty}
< N\delta$. 
Therefore
$$
\eqalign{
\mu_{R\cup \{n\}}(x,a)
\close {N\delta} &
p_{\{n\}}(a)p_R(x) \cr
 \close {N\delta}&  p_{\{n\}}(a)\nu_R(x) \quad
\forall\, R\in \c R, x \in A^R, a \in A.
}
$$

\noindent This means (since we are using $l_{\infty}$-norms) that
for all 
$a \in A$

$$
\{\mu_{R\cup \{n\}}(\cdot,a)\}_{R \in \c R}\close {2N\delta}
p_{\{n\}}(a)\{\nu_R\}_{R \in \c R},
$$

\noindent so using $\|B\| \le C'$ and linearity of $B$ we obtain

$$
\si(\cdot,a) \close {2C'N\delta} p_{\{n\}}(a)
B(\{\nu_R\}_{R \in \c R}) = p_{\{n\}}(a)\pi_{\bar R}\lambda.
\eqno (4)
$$

\noindent Since $\pi_{\bar R}\lambda$ is $\beta$-independent 
and $|{\bar R}|<N$, 
for each
$x \in A^{\bar R}$ we have
$$
\pi_{\bar R}\lambda (x) \close {N\beta}
p_{\bar R}(x) > \alpha^N.
$$
 Taking $\beta =  {\alpha^N  \over 2N}$ it follows that
$\pi_{\bar R}\lambda(x) > \alpha^N/2$ for all $x \in A^{\bar R}$.
Combining this with (4) and taking 
$\delta = {\beta \alpha^N \over 4C'N}$
we get

$$
\left| {\si(x,a) \over \pi_{\bar R}\la(x)} - p_{\{n\}}(a)\right|
    <  {4C'N\delta \over \alpha^N} =  \beta  \qquad
         \forall x \in A^R,\ a \in A.              \eqno   (5)
$$

\noindent Since $p_{\{n\}}(a) \ge \alpha$ and $\beta \le \alpha/2$
this  shows that $\si$ takes only positive values so it is
 a probability measure.  (5) also shows that
$P^{\{n\}} \perp_{\beta,\si}P^{R}$, establishing (3) and concluding
 the proof of  Claim A and Proposition 1.

\bigskip

\noindent {\sl 3. Relatively Independent  Iterates of a Partition}
 
\bigskip

If $T$ is an ergodic invertible measure-preserving tansformation of $X$,
$P$ is a partition of $X$ and
$K$ is any subset of $\Z$ then 
the partition $P^K$, indexed by $A^K$, is
defined by 
$$
[P^K(x)](k) = P(T^k x), \hbox{ for } x \in X, k \in K,
$$
Informally, $P^K(x) = \set{P(T^k x)}_{k \in K}$, the function $P$ evaluated
 along the
$K$-orbit of $x$. $P^K(x)$ is called the $P,K$-name of $x$.
Note that  the atoms of $P^K$ are the atoms of the common refinement of the
partitions $\set{T^{-k}P: k \in K}$. For this reason we it is convenient
 for
us in 
this section to 
restate Theorems 1 and 2 in terms of sequences of
 {\sl negative} times which are
mixing, that is we replace $T$ with $T^{-1}$. Note also that
$\set{n_i}$ is relatively mixing if and only if $\set{-n_i}$ is
relatively mixing.

We will say that $P$ is independent over $K$
if the partitions $\set{T^{-k}P, k \in K}$
 are jointly independent, equivalently
$\dist P^K$ is the product of its one dimensional marginals. We will use the same
terminology in the relative setting, that is with respect to the fiber measures
$\mu_\om$. Similar observations apply to approximate independence,
for example $P^{\set l} $ is $\de$-independent of $P^K$ if and only if
$T^{-l}P$ is $\de$-independent of  $\bigvee_{k \in K}T^{-k}P$.

\proclaim Theorem 1.
Let 
$(X,\c B,\mu,T)$ be an ergodic system with a
factor 
$(\Omega, \c F, \prP ,\theta)$ 
where $\mu$ decomposes as $\mu=\int \mu_\omega d\prP$.
Assume $\{-n_i\}_{i=0}^{\infty}$ is an $\Omega$-mixing sequence for
$T$ and $P$ is a finite 
partition of $X$ for which there is an $\al > 0$ such that
$$
\hbox { for} \prP\hbox{-a.e.} \omega\in \Omega 
\forall p \in P \hbox{ we have }
 \mu_{\omega}(p)>\alpha.
$$ 
Then for every $\ep >0$ there exist a partition $Q$ of $X$,
a set $E \subset \Omega$, $\prP (E) < \ep$,
and 
a subsequence $\{m_i\}$ of $\{n_i\}$, such that
$$
\dist_{\mu_{\omega}} Q = \dist_{\mu_{\omega}} P 
{\hbox { for }} \prP{\hbox{ -a.e. }}\omega, 
\eqno (1)
$$
$$
d_{\mu_{\omega}} (P,Q) < \ep {\hbox {  for }} \prP{\hbox{ -a.e. }} 
\omega
\eqno (2)
$$
and
$$ P
{\hbox { is independent over } \set{m_i}
\hbox{ with respect to }}
\mu_{\omega} {\hbox { for all }} \omega \not \in E.
\eqno (3)
$$

\vskip5mm

\proof
There are two cases: $\theta$ has a set of periodic points of positive
measure, and $\theta$ is 
aperiodic.
In the first case, by the ergodicity of $\theta$ it must be a rotation
on a finite number $k$ of points. 
By dropping to a sub-sequence we may assume that the $n_i$ are all congruent
modulo $k$, say $n_i = n_i' + r$, where the $n_i'$ are multiples of
$k$. The independence we are aiming for over a subsequence $n_{i_j}$
(which is now on {\sl every} fiber) is then
equivalent to independence over $n_{i_j}'$, so we may as well assume that
the $n_i$ are themselves multiples of $k$, say $n_i = q_i k$.
For each $\om$ $T^k$ preserves the measure $\mu_\om$ and
$\set{q_i}$ is mixing for $(T, \mu_\om)$. Fix any $\om_0$ and apply 
the non-relative version of Theorem 1 ([BJ] Theorem 2, also
 [JRW] Theorem 5) to $(T^k,\mu_{\om_0})$ to
modify $P$ by a small amount on the fiber over $\om_0$, obtaining
a  new partition and
a subsequence of $\set{q_i}$ which give the desired independence 
with respect to $\mu_{\om_0}$. Now repeat the process in the fiber over
$\th \om_0$ with the new partition and subsequence to obtain the desired
independence with respect to $\mu_{\th \om_0}$. Iterating the process until
we get to $\th^{k-1}\om_0$ concludes the argument in the periodic case.

In the aperiodic case
we will proceed inductively starting with $Q_0=P$, $m_0=0$,
${\tilde E}_0=\varnothing$ and building the sequences
$\{Q_j\}$, $\{m_j\}$ and $\{E_j\}$ so that for all $j \ge 1$, 
$\prP (E_j) < \ep / 2^j$ and
\medskip
\item{(1$'$)}\quad $\dist_{\mu_{\omega}}Q_j = \dist_{\mu_{\omega}}Q_{j-1} 
{\hbox { for }} \prP{\hbox{-a.e. }}\omega $,
\smallskip
\item{(2$'$)} \quad
$
d_{\mu_{\omega}} (Q_j,Q_{j-1}) < \ep / 2^j {\hbox { for }} 
\prP{\hbox{-a.e. }} 
\omega
$,
\smallskip
\item{(3$'$)}\quad
$Q$
 is jointly independent over $ \set{m_i}_{i=1}^j$ with respect to
$\mu_{\omega}$ for all $ \omega \not \in \tilde E_j$  where
 $ \tilde E_j:=\tilde E_{j-1}\cup E_j 
= \cup_{i=0}^{j} E_i$.
\medskip
Assume for a moment that this construction has been carried out.
Then defining 
$$E = \cup_{j=0}^{\infty} {\tilde E}_j = \cup_{j=0}^{\infty} E_j$$ 
(so $\prP(E) < \ep$) and $Q=\lim Q_j$ we get (1), (2) and (3)
satisfied. The limit is taken in the symmetric difference metric  
in the space of
partitions. The space is complete,
and the sequence of partitions is  Cauchy  by 
($2'$). Therefore we are done.

The proof of the theorem is thereby reduced to the induction step.
For convenience we now get rid of unnecessary subscripts and constants
and restate in equivalent form what remains to be proved.
Use the assumptions of the theorem and let $K$ be  a subset of 
$\N\cup \{0\}$ containing $0$ so that 
 $P$ is independent over $K$ with respect to
$\mu_{\omega}$
for all $\omega$ outside some $E \subset \Omega$.

Given $\ep > 0$ we are looking for $m \in \{n_i\}$, $m>\max  K$, 
 a partition $Q$ of $X$
and $\tilde E \subset \Omega$, $\prP(\tilde E) < \ep$, such that
\medskip
\item{($1''$)} \quad
$
\dist_{\mu_{\omega}} Q = \dist_{\mu_{\omega}}P$ for $ \prP{\hbox{-a.e. }}\omega, 
$
\smallskip
\item{($2''$)} \quad
$
d_{\mu_{\omega}} (P,Q) < \ep {\hbox { for }} \prP{\hbox{-a.e. }} 
\omega,
$
\smallskip
\item{($3''$)} \quad
$
Q{\hbox { is independent over  }} K\cup \set m \hbox{ with respect to }
\mu_{\omega} \ \forall \omega \not \in E\cup \tilde E
$  

\medskip
Denote $|K|=k$ and let 
$\eta =  {1 \over 10} 
{\alpha^{k+1} \over 2}
\de \ep $
where
$\delta = \delta (\alpha, k^2+1)$ from 
Proposition 1 (in section 2) -- this choice 
of $\eta$ will become clear later on. By the hypothesis,
 we know that,
for $\omega \in \Omega
\backslash E$, $ P$ 
is independent over $K$ with respect to
$\mu_{\omega}$, 
in other words
$\dist_{\mu_\om} P^K$ is the product of its one-dimensional
marginals.
Using the mixing property of the sequence $\{-n_i\}$ pick $m \in \set{n_i}$
such that
$m> \max K$ and for  $\omega \notin  E_1$, 
$\prP(E_1)<{\ep \over 10}$,
$$
P^m \perp_{\eta, \mu_{\omega}} P^K.
$$

Let $B_{\Om}
\subset \Om$ be the base of a Rokhlin tower $R$ for $\theta$ of height
$M$ where
$M > 10 m /\ep$, 
with error set $$E_2 = \Om \backslash \bigcup \set{\th^{i}B_\Om: 0 \leq i < M}$$
 of measure less than $ {\ep \over 10}$.
The tower $R$ can be lifted to $X$ by $\pi^{-1}$, where $\pi$ is the
factorization map.
This tower in $X$ can informally be regarded as an ``integral of towers''
with bases 
$B_{\omega}= \pi^{-1}\set{\om},\  \om \in B$, and this observation inspires the
construction that follows.
Let $\bar M= [0,M-1]$.
With $\omega \in B_{\Omega}$ fixed, its $\bar M$-orbit 
$\{\theta^j\omega\}_{j=0}^{M-1}$ gives rise to the $T$-tower
$R_{\omega}$ whose levels are the fibers $\pi^{-1}(\theta^j\omega)$
endowed with measures $\mu_{\theta^j\omega}$.

Let $\nu_{\omega}= \dist _{\mu_{\omega}} P^{\bar M}$, a probability
measure on $A^{\bar M}$.
Denote $K'=K \cup \{m\}$. For all 
$j$ that satisfy $j+K' \subset \bar M$
let $\nu_{\omega,j} =\pi_{j+K'}\nu_{\omega}$. 
For those values of $j$ for which 
$\theta^j \omega \not \in E \cup E_1$ the measure
$\dist \nu_{\om,j}$ on $A^{j+K'}$
is $\eta$-independent,  by our choice of $m$ and because
 its marginal on $A^{j+K}$ is exactly independent.
We now aim at perturbing $P$ to obtain $Q$ that satisfies
$(1'')$ and $(2'')$ and for which
$\dist_{\mu_\om}Q^{j+K'}$ is exactly independent  for
all combinations of $\omega$ and $j$ corresponding to 
$\theta^j \omega$ outside of a set slightly larger than $E \cup E_1$.

Fix an $\omega \in B_{\Omega}$ and consider the tower $R_\omega$.
Denote its base $\pi^{-1}\set{\om}$ by $B$.
Partition $B$ into two parts $B=B_0 \cup B_1$, 
$\mu_\omega (B_0) = {\ep \over 10}$,
independently of $\bar M$-names of the points in $B$ (i.e. after proper 
normalization
$\dist_{\mu_\omega } P^{\bar M}$ on each of $B_0$ and $B_1$
is identical to that on $B$).   
This can be done since the fibers of $\pi$ are non-atomic by the mixing
property of $T$.

We will construct the required partition $Q$ by
re-assigning
points on the $M$-tower above $B_0$ to the letters of the alphabet
$A$ (the technique known as ``painting names on towers''), while 
retaining this assignment on the $M$-tower above 
$B_1$. We now explain how the re-assignment is done. 
(Recall that $\omega$ is fixed
for the time being.)

Let $m_{\omega,j}$ denote the projection of $\nu_\om$ onto
$A^{\set{j}}$, a probability measure on $A^{\set{j}}$ which is the same as 
$\dist_{\mu_{\th^j \om}}P$
 (up to identification of
 $A^{\set{j}}$ with
 $A^{\set{0}}$).
Of course
the one-dimensional marginal of $\nu_{\omega,j}$ on $A^{\{j+k\}}$
for any particular $k$ is $m_{\omega, j+k}$.
For all $j$ satisfying
 $$ 0\le j\le M-m \hbox{ and }
  \th^j \om  \notin E \cup E_1  \eqno (*)$$
 we 
define a ``correcting measure'' $ \xi_{\omega,j}$ on $A^{j+K'}$  by  the 
following equation:  
$$
(1-({\ep}/{10})) \nu_{\omega,j} + ({\ep}/{10}) \xi_{\omega,j}
=\prod_{k \in K'} m_{\omega,j+k}.
$$
The idea is that blending a small fraction of 
$\xi_{\om,j}$ into $\nu_{\om,j}$ corrects it  from approximate 
to exact independence.
 A priori  $\xi_{\omega,j}$ is a signed measure
but we shall see shortly that it is in fact positive. Observe that
for any $k \in K'$
the marginal of   $\xi_{\omega,j}$ on $A^{j+k}$ is 
$m_{\om, j+k}$, since that is the case for both $\nu_{\om,j}$ 
and $\prod_{k \in K'} m_{\omega,j+k}$. 


The definition of $\xi_{\om,j}$ is equivalent to
$$
\xi_{\omega,j}-\prod_{k \in K'}  m_{\omega,j+k}
=(({10}/{\ep}) -1) 
\bigl(\prod_{k \in K'}  m_{\omega,j+k} - 
\nu_{\omega,j}\bigr).
\eqno (4)
$$
By our choice of $m$ the $l_{\infty}$-norm of the right hand
side of (4) is less than 
$$
{10 \over \ep}\eta =  \de  {\al^{k + 1} \over 2 } <  {\al^{k + 1} \over 2 }.
$$ 
Since the minimal value of $\prod_{k \in K'}  m_{\omega,j+k} $ 
 on  singletons in $A^{j+K'}$ is at least $\al^{k + 1}$ we conclude that
$\xi_{\om ,j}$ is indeed positive.

Projecting the relation between the measures from $A^{j+K'}$
to $A^{j+K}$ we conclude that
$\pi_{A^{j+K}}\xi_{\omega,j}$ is the product of its
one-dimensional marginals \break $ m_{\omega,j+k}$,
since this is true for $\pi_{A^{j+K}}\nu_{\omega,j}$.  
In addition, we claim 
that for  all $j$ satisfying (*) the measure
$\xi_{\omega,j}$
is $\delta$-independent.
In view of the
 mutual independence of the first $k$ 
(out of the total of $k+1$)
one-dimensional marginals of this measure it would be enough to 
prove that for all
$ x \in {A^{j+K}}$ and $ a \in A$ we have
$$
\bigl|{\xi_{\omega,j}}(x,a)/\pi_{A^{j+K}}\xi_{\omega,j}(x)-
\pi_{A^{\{j+m\}}}\xi_{\omega,j}(a)
\bigr| < \delta.
$$
(As in Section 2 we view $y\in {A^{j+K'}}$ as a pair $y=(x,a)$,
$x \in {A^{j+K}}, a \in {A^{\{j+m\}}}$.)
This is achieved by our choice of 
$\eta$,
as the following calculation
shows. Evaluate (4) at $(x,a)$ and then divide by 
$$\pi_{A^{j+K}}\xi_{\omega,j}(x)
=\pi_{A^{j+K}}\nu_{\omega,j}(x)
=\prod_{k \in K} m_{\omega, j+k}$$
 to obtain
$$
\eqalign{
&\bigl|{\xi_{\omega,j}}(x,a)/\pi_{A^{j+K}}\xi_{\omega,j}(x)-
\pi_{A^{\{j+m\}}}\xi_{\omega,j}(a)
\bigr|     \cr
&=\bigl( {10 \over \ep}-1\bigr)
\Bigl| {\nu_{\omega,j}(x,a) \over \pi_{A^{j+K}}\nu_{\omega,j}(x)}-
\pi_{A^{\{j+m\}}}\nu_{\omega,j}(a)
\Bigr| \cr
&<{10 \over \ep}\cdot{1 \over 10}
\delta\ep {\al^{k+1} \over 2} \leq \delta.
}
$$

 Proposition 1 of Section 2 can now be applied to the 
family $\c F$ of those shifts $j+K'$ of $K'$ for which
$j$ satisfies (*),
endowed by the measures 
${\xi_{\omega,j}}$. 
({\sl F}) holds for $N=k^2+1$. (M1), the consistency condition for the measures
${\xi_{\omega,j}}$, is implied by the consistency of 
${\nu_{\omega,j}}$. (M2) holds because the 1-dimensional marginals of 
${\xi_{\omega,j}}$ are the same as for ${\nu_{\omega,j}}$. (M3) has just
been demonstrated. So there is a measure $\xi_{\omega}$ 
on $A^{\bar M}$ that is a common extension of all the measures
in the family. Write the names on the levels of the $T$-tower 
above $B_0$ as prescribed by $\xi_{\omega}$. Explicitly what this means is
we
choose a partition $R$ of $B_0$ indexed by $A^{\bar M}$
so that $\dist (R|B)=\xi_{\omega}$ and thenwe
 partition the  tower over $B_0$ according to the rule:
for each $x \in B_0$ and $j \in \bar M$ let 
$Q(T^j x)=(R(x))(j)$. 
To complete the construction of the partition $Q$ on the whole 
$M$-tower above $B$ recall that this partition coincides with $P$
above $B_1$.
   
The procedure we described depends on the point $\omega \in
B_{\Omega}$. Performing it for every
$\omega \in B_{\Omega}$ we get a partition $Q$ of the whole
$\pi^{-1}(\Omega \backslash E_2)$.
On $\pi^{-1}(E_2)$ just declare $Q:=P$.
One only needs to assure that $Q$ can be chosen $\mu$-measurable. This 
follows from the fact that for any measurable $A \subset X$ and  $0 <t <1$ there
is a set $B \subset X$ such that $\mu_\om B =t$ for a.a. $\om$ and
$A$ and $B$ are independent with respect to $\mu_{\om}$ for a.a. $\om$.
 Using this and following every step of the procedure it is
easy to convince oneself that everything can be done in a $\mu$-measurable
way. 

Having done all this, we obtain a partition $Q$ of $X$ that clearly
satisfies $(1'')$ and $(2'')$. $(1'')$ holds
because the 1-dimensional marginals
of the measures
$\nu_{\om,j}$ and $\xi_{\om,j}$ are the same.
$(2'')$ is satisfied even with $\ep/10$ instead of $\ep$.
By virtue of the construction, $(3'')$ is satisfied for all $\omega$
outside $E\cup E_1 \cup E_2 \cup E_3$.
$E_1$ and $E_2$ have been defined above, while $E_3$ is the set of
$m$ uppermost levels in the $M$-tower for $\theta$  with
base $B_{\Omega}$. So $\prP (E_3) < m/M < \ep /10$.
Therefore taking $\tilde E =E_1 \cup E_2 \cup E_3$ we
obtain the required version of $(3'')$.

\endproof

\noindent {\sl Remark.}\   The proof shows that,
  as in  
[JRW]and [JB],
 the independence can be 
achieved not
just along the sequence $\set{m_i}$, but also along the
IP-set which it generates, namely the set of all finite sums of the form
$\sum_{i \in F}m_i$, where $F$ is a finite subset of $\N$.

\proclaim Theorem 2.
With the same hypotheses as in theorem 1
  there is 
a $Q$  and a subsequence $\{m_i\}$ such that
$\dist_{\mu_\om} Q = \dist_{\mu_\om} P$ for a.a. $\om$,
$d(P,Q) < \ep$ and 
$Q$ is independent over $\set{m_i}$ with respect to $\mu_\om$ for
almost every $\om$.

The proof is a modification of the proof of Theorem 1, using the following lemma.

\proclaim Lemma 3.
Suppose $K$ is a finite subset of $\Bbb N$ and $P$ is a partition
which is independent over $K$ with
respect to $\mu_\om$ for all $\om \notin E $ where $E$ is some subset  of $\Om$.
Then for any $\ep >0$ there is a partition $Q$ such that $P$ and $Q$ agree
 on all fibers not in $KE=\cup\set{T^k E: k \in K}$,
$\dist_{\mu_\om}P = \dist_{\mu_\om} Q$ for all $\om$ and
there is a subset
$F \subset \Om$ with $\mu(F) <   \ep$ such that
$Q$ is independent over $K$ with respect to
$\mu_\om$ for all $\om \notin F$.

\proof
Let $B \subset \Om$ be the base of a Rohlin tower of height $n \gg k=\max K$.
 Fix $\om_0 \in B$.
As in the proof of Theorem 1 we will describe how to change
$P$ on the tower over $\pi^{-1}\set{\om_0}$, with the understanding
that this should be done simultaneously for  each 
$\om_0 \in B$, in a measurable varying way.
 
Fix any $i_1$ such that  $0 \leq i_1 \leq  n-k$ and
$\theta^{i_1}\om_0 \in E$.
 Now change $P$ on each $\pi^{-1}(\theta^{i_1 +k}\om_0), \ k \in K$,
retaining the same distribution on each of these fibers, to obtain
a new partition $P_1$ so that, with respect to $\mu_{\om_0}$,
$P_1$ is independent over $ i_1 + K$ and in addition
 $P_1^{i_1 +K}$ is independent of $P_1^{[0,n] \backslash (i+K)}$.
This is possible simply because the fibers are  non-atomic.

 Of course we will 
then have  that  $ P_1$ is independent over $K$ with respect to
$\mu_{\th^{i_1} \om_0}$. In addition,
 for any $i \neq i_0$  such that  $0 \leq i \leq  n-k$
 and
$\theta^{i}\om_0 \notin E$, we will still have that
 $ P_1^K$ is  independent over $K$ with respect to
$\mu_{\th^i \om_0}$.
 Indeed, with respect to  $\mu_{\om_0}$, $P_1$ is independent
over $S_1 = i+K \cap (i_1+ K)^c$, since $P$ has this property while $P_1$ agrees with
$P$ over the relevant fibers. Also with respect to $\mu_{\om_0}$,
$P_1$ is independent over  $ S_2 = i+K \cap (i_1 +K)$, since
$P_1$ is actually independent over all of $i_1 + K$ by construction. Finally
$P_1^{S_1}$ and $P_1^{S_2}$ are independent with respect to
$\mu_{\om_0}$, again by construction of $P_1$.
This shows that $P_1$ is independent over $i +K$ with respect to
$\mu_{\om_0}$ which implies that $P_1$ is independent over $K$ with respect
to $\mu_{\th^i \om_0}$

 Now replace $i_1$  with any $i_2, 0 \leq i_2 \leq n-k$ 
such that we still do not have the desired independence 
with respect to $\mu_{\th^{i_2} \om_0}$
and modify $P_1$ to $P_2$
as before to achieve the independence
with respect to $\mu_{\th^{i_2}\om_0}$.
Continuing in this way we will arrive at a partition $Q$ which has the desired
independence with respect to $\mu_{\th^i \om_0}$ for every
$i$, $ 0 \leq i < n-k$.
 Note also that $Q$ differs from
$P$ only on fibers over over $KE$.
Doing this for all $\om_0 \in B$ we obtain the desired independence
with respect to $\mu_{\om}$ for all $\om \in \th^i B$, $0 \leq i < n-k$ and
the measure of this set of $\om$'s can be made as small as we please.

\endproof

\proclaim Corollary 1. Lemma 3 holds with $\ep=0$.

\proof Iterating Lemma 3 with a summable sequence of $\ep$'s produces a sequence
of partition which, by the 
Borel-Cantelli lemma, eventually stabilizes over
 $\om$ for a.a. $\om$. This gives
a limiting partition $Q$ with the desired properties.

\endproof

Now it is to prove Theorem 2:
in the proof of theorem 1 at each stage, after achieving the
 independence of $P_k$ over $\set{m_1, \ldots, m_k}$ with respect to
$\mu_\om$
for a large set of $\om$,
use Corollary 1 to change $P_k$  over a small set of $\om$
to achieve the desired independence over a.a. $\om$. This concludes the proof
of Theorem 2.

                       \endproof

\bigskip

We now sketch how to remove the ergodicity
assumption in theorems 1 and 2.
Let us first observe that it is enough to prove the result in the
case when $\th$ is periodic with period $n \in \set{1,2, \ldots, \infty}$.
Indeed, supposing that we have done this,
 partition $\Om$ into the  invariant  sets $\Om_n$,
$n=1,2, \ldots,  \infty$
where $\th$ has period exactly  $n$. We start in $\Om_{\infty}$ and find a
 $Q_{\infty}$ uniformly close to $P$ over $\Om_{\infty}$ and equal to
$P$ over $\Om \backslash \Om_{\infty}$,
and a first
subsequence $\set {m_i}$ which gives the desired independence
in the fibers over most of $\Om_\infty$. Then we work over
$\Om_1$ and find $Q_1$ close to
$Q_{\infty}$ and a further 
subsequence which gives the desired independence on most of
$\Om_1$. We continue in this way, perturbing the partition and refining
the subsequence, for finitely many steps until we have
exhausted all but a small invariant set $E \subset \Om$ consisisting of some
tail of the sequence $\set{\Om_n}$.
For Theorem 1 we only need
independence on all but a small subset of $\Om$ so we are  done.
In the case of theorem 2,
if we are willing to relax 
 the requirement that the perturbed partition  
have the same distribution as the original one on all fibers,
 we can get the desired independence on fibers
over $E$ by simply making the partition trivial on those fibers.

So, we may now assume that $\th$ has constant period $k$, possibly infinite.
In both theorems 1 and 2 the case $k = \infty$ is identical to the ergodic case
since aperiodicity is all that is needed for the Rohlin lemma.
In the case when $k$ is finite, Theorem 1, the same 
congruence trick used in the ergodic case
reduces us to working with $T^k$ so we may assume that $\th = {\rm id}_{\Om}$. 
So we have reduced ourselves to proving  the non-ergodic versions of
Theorems 1 and 2 in the case when $\th$ is the
identity. In the case of Theorem 1 this is just a slightly souped up
version of the absolute theorem where instead of having just one
weakly mixing system one has a whole measurable field of them.
It is quite straightforward to prove this by suitably modifying the
proof of the absolute result.  To actually carry it out  would require us
to
 delve  into the proof
of the absolute result. This is the reason that we chose to prove theorems
1 and 2 formally only in the ergodic case.
In the case of Theorem 2 we again just make the partition trivial in
the small set of fibers where Theorem 1 does not give independence.

The observant reader will have noticed that Theorems 1 and 2 do not
quite generalize Krengel's theorem in that we don't insist that the
subsequence $\set{m_i}$  start with $m_1 = 0$. The place where our
argument would break down if we wanted $m_1=0$ is at the congruence
trick used to deal with the periodic case. The results may well be true
with $m_1= 0$ and in any case our argument shows that they do hold
at least when $\th$ is aperiodic.

\bigskip

 {\sl Acknowledgments.}
We would like to thank
B. Weiss and Y.Kifer for helpful discussions that led to
improvement of the exposition.

\bigskip

\centerline{References}
\bigskip
\item{ [BJ]}
B. Begun and A. del Junco,
{\sl Amenable groups, stationary measures and
partitions with independent iterates}
to appear in {\sl Israel J. Math.}
\smallskip
\item{ [F]}
H. Furstenberg,
{\sl Recurrence in Ergodic Theory and Combinatorial Number Theory},
Princeton University Press, 1981.

\smallskip
\item{[G]}
E. Glasner,
{\sl Ergodic Theory via Joinings\/},
Math. Surveys and Monographs {\bf 101},
Amer. Math. Soc. 2003.
\smallskip
\item{ [JRW]}
A. del Junco, K. Reinhold, and B.Weiss,
{\sl Partitions with independent iterates along IP-sets},
{ Ergodic Theory Dynam. Systems\/} {\bf 19} (1999), 447-473.
\smallskip
\item{ [K]}
U. Krengel,
{\sl Weakly wandering vectors and weakly independent partitions},
{ Trans. Amer. Math. Soc.\/} {\bf 164} (1972), 199-226.
\smallskip
\item{ [KW]}
Y. Kifer and B. Weiss,
{\sl Generating partitions for random transformations},
{ Ergodic Theory Dynam. Systems\/} {\bf 22} (2002), 1813-1830.
\smallskip
\item{ [R]}
D. J. Rudolph,
{\sl Pointwise and $L^1$ mixing relative to a sub-sigma algebra},
{ Illinois J. Math.\/} {\bf 48} (2004), 505-517.

\end

%% file: macros.tex
\def\dist{{\,{\rm dist}\,}}

\def\proof{{\noindent{\bf Proof:} \quad }}
\def\endproof{\hfill{\vrule height4pt width6pt depth2pt} \vskip .2in}

\def\to{\rightarrow}

\def\vec3#1#2#3{(#1,#2,#3)}

\def\close#1{\ {\mathop{\sim}\limits^{#1}}\ }

\def\c#1{{\cal #1}}

\def\dist{{\,{\rm dist}\,}}

\def\proof{{\noindent{\bf Proof:} \qquad }}

\def\set#1{\{#1\}}

\def\al{\alpha}

\def\de{\delta}
\def\ep{\epsilon}

\def\th{\theta}

\def\la{\lambda}
\def\si{\sigma}

\def\ph{\varphi}

\def\om{\omega}

\def\bar#1{\overline #1}

\def\today{{\ifcase\month\or January\or February\or March \or April
\or May\or June\or July \or August \or September \or
October \or November \or December\fi}\  \number\day,\
\number\year}


\def\su #1 #2 #3 {#1_{#2}^{#3}}

\def\abs#1{|#1|}

\def\X{(X, \c B, \mu,T)}

\def\N{\Bbb N}
\def\R{\Bbb R}
\def\Z{\Bbb Z}

\def\F{\Bbb F}


